%% THIS IS VOLUME THIRTEEN
%
%
\magnification=1200
%\nopagenumbers
\hfuzz 8pt
\vsize 23truecm
\hsize 16.2truecm
\parindent 15pt
\vbox{\vskip 3truecm}
%%%%%%%%%%%%%%%%%%%%%%%%%%%%%%%%%%%%%%%%%%%%%%%%%%%%%%%%%%%%%%%%%%%%%%%%%
\font\ams=msbm10
\font\sym=msbm8

%%%%%%%%%%%%%%%%%%%%%%%%%%%%%%%%%%%%%%%%%%%%%%%%%%%%%%%%%%%%%%%%%%%%%%%%%%
\def\lx{\hbox{\sym\char'130}}
\def\ln{\hbox{\sym\char'114}}
\def\lk{\hbox{\sym\char'113}}
\def\lt{\hbox{\sym\char'124}}
\def\lh{\hbox{\sym\char'110}}
%%%%%%%%%%%%%%%%%%%%%%%%%%%%%%%%
\def\cline{\hbox{\ams\char'103}}

\def\hline{\hbox{\ams\char'110}}
\def\kline{\hbox{\ams\char'113}}
\def\lline{\hbox{\ams\char'114}}
\def\mline{\hbox{\ams\char'115}}
\def\nline{\hbox{\ams\char'116}}
\def\tline{\hbox{\ams\char'124}}
\def\xline{\hbox{\ams\char'130}}
\def\sksix{\noalign{\vskip 6pt}}
\def\skthr{\noalign{\vskip 3pt}}
\def\({[}
\def\){]}
\def\bk{\hfill\break}

\font\ninerm = cmr9
%%%%%
% Changed definitions
%%%%%
\def\rar{\Bigl\rangle}
\def\lal{\Bigl\langle}
\def\blacksq{\vrule height 4pt width 5pt depth 2pt}
%%%%%%%%%%%%%%%%%%%%%%%%%%%%%%%%%%%%%%%%%%%%%%%%%%%%%%%%%%%%%%%%%
\newdimen\offdimen
\def\offset#1#2{\offdimen #1
   \noindent \hangindent \offdimen
   \hbox to \offdimen{#2\hfil}\ignorespaces}
%%%%%%%%%%%%%%%%%%%%%%%%%%%%%%%%%%%%%%%%%%%%%%%%%%%%%%%%%%%%%%%%%
%
\baselineskip 10pt
\centerline{\bf AN ABSTRACT INTERPOLATION PROBLEM AND THE}
\medskip
\centerline{\bf EXTENSION THEORY OF ISOMETRIC OPERATORS%
\footnote{*}{\ninerm Translated from ``Operators in Function
Spaces and Problems in Function Theory",
\bk\vskip -9pt
pp. 83--96. (Naukova Dumka, Kiev, 1987, Edited by V. A.
Marchenko.}}\bigskip\baselineskip 13pt
\centerline{\bf V. Katsnelson, A. Kheifets and
P. Yuditskii}\vskip 1truecm
\baselineskip 11pt
{\leftskip 35pt\rightskip 35pt
{\ninerm
The algebraic structure of V.P. Potapov's Fundamental Matrix
Inequality (FMI) is discussed and its interpolation meaning
is analyzed. Functional model spaces are involved.  A
general Abstract Interpolation Problem is formulated which seems
to cover all the classical and recent problems in the field and
the solution set of this problem is described using the
Arov--Grossman formula. The extension theory of isometric
operators is the proper language for treating interpolation
problems of this type.}\par}
\vskip 1truecm\baselineskip 15pt
\noindent
{\bf 1. INTERPOLATION DATA AND EXAMPLES}
\medskip
Let $\lline$ and $\lline^{\prime}$ be Hilbert spaces. We shall
denote by $B(\lline,\lline^{\prime})$ the class of operator-%
valued functions which are holomorphic in the unit disk
$\vert\zeta\vert<1$ and whose values are contractive operators
from $\lline$ into $\lline^{\prime}$.  Let $\xline$ be a linear
space. We do not suppose that $\xline$ is  endowed with a
topological structure. Let  $D$  be a sesquilinear form in
$\xline$ and let $T$ be a linear operator on $\xline$. Let $E$
and $M$ be linear operators from $\xline$ into $\lline$ and
$\lline^{\prime}$, respectively.
\par\noindent
We assume that the operators and the sesquilinear form
are linked through the so-called
{\sl Fundamental Identity }  (F\thinspace I):
$$ D(x,y)-D(Tx,Ty)=\langle Ex,Ey\rangle _{\ln}-
\langle Mx,My \rangle _{\ln^{\prime}}\ .\eqno (FI) $$
The Fundamental  Identity  must be fulfilled for arbitrary $x$
and $y$ in $\xline$. This framework arose from the study of a
number of interpolation problems. We shall give some examples.
\par
EXAMPLE 1.~~(The Nevanlinna-Pick Problem).
\par\noindent
{\sl Interpolation Data}: {A sequence $\{\zeta_k\}_{1\leq k\leq
\infty}$ of complex numbers $(\vert\zeta_k\vert<1)$ and a
sequence of contractive operators $\{s_k\}_{1\leq k\leq\infty}$,
acting from $\lline$ into $\lline^{\prime}$.}\par\noindent
{\sl A solution} \ $s(\zeta)$ of this problem is {an arbitrary
holomorphic function from the class $B(\lline,\lline^{\prime})$
which satisfies the
interpolation conditions}\ \ $s(\zeta_k)=s_k$.~
{\sl It is required} {to give a criteria for the existence of
solutions and to describe the set of all solutions of the
interpolation problem.}
\par
In this example, $\xline$ is the space of all infinite sequences
whose entries are vectors from $\lline$ such that only a finite
number of them do not vanish:
$$ x=\{\ell_1,\ell_2,\cdots,\ell_n,0,\cdots,0,\cdots \} \ . $$
The operators $T,\ E,\ M$ and the sesquilinear form $D$ are
defined by the formulas
$$ Tx =\{\zeta_1\ell_1,\zeta_2\ell_2,\cdots,\zeta_n\ell_n,0,
\cdots,0,\cdots\} \ ,$$
$$ Ex=\sum_{1\leq k<\infty}\ell_k~,~~~~Mx=\sum_{1\leq k<\infty}
s_k\ell_k \ , $$
$$ D(x,x)=\sum_{1\leq j,k<\infty}\lal{I-s^*_js_k\over 1-
\overline\zeta_j\zeta_k} \ \ell_k, \ell_j\rar_{\ln}\ .  $$
\par\noindent
The Fundamental Identity   can be verified directly. The
condition $D(x,x)\geq 0\ (\forall\ x\in\xline)$ is necessary and
sufficient for the solvability of this problem.
\par
{EXAMPLE 2~(The Sarason Problem).}~Let $\theta$ be an inner
function, let $\kline_\theta:=H^2\ominus\theta H^2$, and let
$P_\theta$ be the orthoprojection from $H^2$ onto $\kline_
\theta$.
\par\noindent
{\sl The interpolation data} for this problem are the operator
$T=P_\theta t_{\vert{\lk}_\theta}$ i.e., the compressed shift on
 $\kline_\theta$, and a contractive operator $W$ which
commutes with $T$.
\par\noindent
{\sl A solution} of this problem is an arbitrary holomorphic
function $w(\zeta)$ from the class $B(\lline,\lline^{\prime})$
which satisfies the {\sl interpolation condition} $W=P_\theta w
\vert_{\kline_\theta}$.
{\sl It is required} {to describe the set of all solutions to
this problem. The existence of solutions follows from the
contractivity of $W$.}
\par
In this example $\xline$ coincides with $\kline_\theta$, the
operator $T$ is defined in the statement of the problem,
$\lline=\lline^{\prime}=\cline$,
$$ Ex \ {\buildrel def\over{=\hskip -4.0pt =}}\ \langle x,e_*
\rangle\ ,\ \ Mx\ {\buildrel def\over{=\hskip -4.0pt =}}\
\langle Wx,e_*\rangle \ $$
and
$$ D(x,x)\ {\buildrel def\over{=\hskip -4.0pt =}}  \ \langle
(I-W^*W)x,x\rangle \ \ , $$
where $\langle\ ,\ \rangle$ is the inner product in $\kline_
\theta$, the vector $e_*\in\kline_\theta$ is defined by  the
formula $e_*=\(\theta(t)-\theta(0)\)/t$. The Fundamental Identity
is a consequence of the fact that
$$ (I-T^*T)x= e_*\langle x,e_*\rangle $$
for any $x\in\kline_\theta$.
%%%%%%%%%%%%%%%%%%%%%%%%%%%%%%%%%%%%%%%%%%%%%%%%%%%%%%%%%%%%%%%%%
\medskip
\noindent\hangindent 20pt
{\bf 2.}~~{\bf V.P. POTAPOV'S FUNDAMENTAL MATRIX INEQUALITY\bk
AND ITS TRANSFORMATION}\medskip
We shall say that a function  $s\in B(\lline,\lline^{\prime})$
satisfies the V.P. Potapov's Fundamental Matrix Inequality (the
FMI) if for arbitrary $x\in\xline\ ,\ \ell\in\lline$ and
$\vert\zeta\vert<1$,
$$ \left\(\matrix{ D((I-T\overline\zeta)x,(I-T\overline\zeta)x)
&\quad
& \langle \ell,(E-s^*(\zeta)M)x\rangle_{\ln}
\cr \ & \ & \ \cr
\langle (E-s^*(\zeta)M)x,\ell\rangle_{\ln}
&\quad
& \langle {I_{\ln}-s^*(\zeta)s(\zeta)\over 1
-\overline\zeta\zeta}\ \ell,\ell \rangle_{\ln} \cr}\right\)\geq
0\ . \ \eqno(FMI) $$
\par\noindent
This inequality can be rewritten in another (equivalent) form:
$$ \left\( \matrix {D((\zeta I-T)x,(\zeta I-T)x)& \quad
&\langle \ell^{\prime},(s(\zeta)E-M)x,\ell^{\prime} \rangle_
%{\lline^{\prime}}
{\ln^{\prime}}\cr \ & \ & \ \cr
\langle (s(\zeta)E-M)x,\ell^{\prime}\rangle_{\ln^{\prime}}&
\quad & \langle {I_{\ln^{\prime}}-s(\zeta)s^*(\zeta)
\over 1-\zeta\overline\zeta}\ \ell^{\prime},\ell^{\prime}
\rangle_{\ln^{\prime}} \cr}
\right\)\geq 0 \ \eqno(FMI^{\prime}) $$
for arbitrary $x\in\xline\ ,\ \ell^{\prime}\in\lline^{\prime}$
and $\ \vert\zeta\vert<1$.\par
It will be shown, that the equivalence of the FMI and the
FMI$^{\prime}$ is a consequence of the Fundamental Identity  .
\par
In this work we shall use several facts which we
formulate here as propositions 1, 2 and 3 without proofs.
\vskip 0.3 truecm \par
{PROPOSITION 1.}~~(The block-matrix lemma).~
{\sl Let $\hline$ be a Hilbert space and let $A$ be a
selfadjoint positive semidefinite ($A\geq 0$) operator
acting on $\hline$. Suppose that $h_0\in\hline$ and that there
exists a constant $C\geq 0$ such that
$$\vert\langle h_0,h\rangle\vert^2\leq C \|\sqrt A h \| ^2 \c
\qquad \forall\ h\in D_{\sqrt A}\ \ \c\ $$
where $D_{\sqrt A}$ is the domain of the selfadjoint operator
$\sqrt A$. Then there exists a unique vector $g_0\in (Ker A)^
\perp\cap D_{\sqrt A}$ such that $\sqrt A g_0=h_0 .$ Moreover,
$$\Vert g_0\Vert^2\leq C .$$}
We shall use the following notations :
$$ g_0=A^{\(-{1\over 2}\)}h_0 {\rm \ \ \ and\ \ \ }
\langle g_0 , g_0\rangle\ =  \langle A^{\(-1\)} h_0,h_0\rangle .
$$
\vskip 0.3 truecm
{REMARK.}~~The converse statement is obvious : If $A\geq 0\ ,
\ g_0\in D_{\sqrt A} {\rm\ \ and\ \ }$
$$\Vert g_0\Vert^2\leq C , $$
then
$$\vert\langle \sqrt A g_0,h\rangle\vert^2\leq C \|\sqrt A
h \| ^2 $$
for arbitrary $h\in D_{\sqrt A}\ \ .$
\vskip 0.4 truecm\par
{PROPOSITION 2.}~~{\sl Let $s$ be a contractive operator from
$\lline$ into $\lline^{\prime}$, then
$$ s(I_{\ln}-s^*s)^{\(-{1\over 2}\)}=(I_{\ln^{\prime}}-ss^*)^
{\(-{1\over 2}\)}s $$
and
$$ (I_{\ln}-s^*s)^{\(-{1\over 2}\)}s^*=s^*(I_{\ln^{\prime}}-ss^*)
^{\(-{1\over 2}\)}\ . $$}
\noindent
In particular, Proposition 2 contains the assertion that the
domains of the operators on the right and left hand sides of
these equalities coincide.\par
{PROPOSITION 3.}~~{\sl Let $s$ be a contractive operator from
$\lline$ into $\lline^{\prime}$. Then the following statements
are equivalent:
$$ \ell^{\prime}\oplus\ell\quad {belongs\ to\ the\
domain\ of\  the\ operator\ } \left\( \matrix
{I_{\ln^{\prime}} & s \cr s^* & I_{\ln} \cr}\right\)^{\(-{1\over
2}\)} \leqno (1) $$
\vskip -27pt
$$ \ell-s^*\ell^{\prime}\quad {belongs\ to\ the\ domain\
of\  the\ operator\ } (I_{\lline}-s^*s)^{\(-{1\over 2}\)}
\leqno (2) $$
\vskip -27pt
$$ \ell^{\prime}-s\ell\quad {belongs\ to\ the\ domain\
of\  the\ operator\ } (I_{\ln^{\prime}}-ss^*)^{\(-{1\over 2}\)}
. \leqno (3) $$
Moreover, if these properties are in force, then
$$ \eqalign{&
\langle \left\(\matrix{I_{\ln^{\prime}} & s \cr s^* & I_{\ln}\cr}
\right\)^{\(-1\)}(\ell^{\prime}\oplus\ell),(\ell^{\prime}\oplus
\ell)\rangle_{\ln^{\prime}\oplus\ln}  \cr
& =\ \langle \ell^{\prime},\ell^{\prime}\rangle _{\ln^{\prime}}+
\langle (I_{\ln}-s^*s)^{\(-1\)}(\ell-s^*
\ell^{\prime}),(\ell-s^*\ell^{\prime})\rangle _{\ln}\  \cr
& =\ \langle \ell,\ell\rangle _{\ln}+\langle (I_{\ln^{\prime}}
-ss^*)^{\(-1\)}(\ell^{\prime}-s\ell),
(\ell^{\prime}-s\ell)\rangle _{\ln^{\prime}} \ . \cr} $$}
\par\noindent
The latter formulas inply the following:\par
{COROLLARY.} {\sl The inequalities
$$ \langle \left\(\matrix{I_{\ln^{\prime}} & s \cr s^* & I_{\ln}
\cr}\right\)^{\(-1\)}(\ell^{\prime}\oplus\ell),(\ell^{\prime}
\oplus\ell)\rangle_{\ln^{\prime}\oplus\ln} \geq \ \langle \ell,
\ell\rangle _{\ln} $$
and
$$\langle \left\(\matrix{I_{\ln^{\prime}} & s \cr s^* & I_{\ln}
\cr}\right\)^{\(-1\)}(\ell^{\prime}\oplus\ell),(\ell^{\prime}
\oplus\ell)\rangle_{\ln^{\prime}\oplus\ln} \geq \ \langle \ell^
{\prime},\ell^{\prime}\rangle _{\ln^{\prime}} $$
hold.}

\noindent
We are now ready to prove the equivalence of the FMI and the
FMI$^{\prime}$.\par
{PROPOSITION 4.} \ \ {\sl The FMI holds if and only if the
FMI$^{\ \prime}$ holds.}\par
{PROOF.} \ \ In view of the block-matrix lemma, the FMI is
equivalent to the following inequality
$$ \eqalign{&  D((I-T\overline\zeta)x,(I-T\overline\zeta)x)\cr
& -\langle\left\({I_{\ln^{\prime}}-s^*(\zeta)s(\zeta)\over 1-
\overline\zeta\zeta}\right\)^{\(-1\)}
(E-s^*(\zeta)M)x,(E-s^*(\zeta)M)x \rangle \geq 0 \ . \cr} $$
Using Proposition 3 we get
$$ \eqalign{&  D((I-T\overline\zeta)x,(I-T\overline\zeta)x)+
(1-\zeta\overline\zeta)\{\langle Mx,
Mx \rangle -\langle Ex,Ex \rangle \} \cr
& -\langle \left\({I_{\ln}-s(\zeta)s^*(\zeta)\over 1-\zeta
\overline\zeta}\right\)^{\(-1\)}
(s(\zeta)E-M)x,(s(\zeta)E-M)x \rangle\geq 0 \ . \cr} $$
According to the Fundamental Identity  \
$$ \langle Mx,Mx\rangle -\langle Ex,Ex\rangle \ =D(Tx,Tx)-D(x,x)
\ . $$
Combining similar terms in the sesquilinear form $D$ we obtain
the inequality
$$ D((\zeta I-T)x,(\zeta I-T)x)-\langle\left\({I_{\ln^{\prime}}-
s(\zeta)s^*(\zeta)\over 1-\zeta\overline\zeta}\right\)^
{\(-1\)} (s(\zeta)E-M)x,(s(\zeta)E-M)x \rangle \geq 0 \ , $$
which is equivalent to the FMI because of the block-matrix lemma.
Hence Proposition 4 is proved.\hfill\blacksq\par
Our aim now is to extract interpolation information from the FMI.
We shall apply the method used earlier by one of the authors in
\(3\).\par
We begin with an illustrative example: Let $\zeta_0$ be an
eigenvalue of $T$, let $x_0$ be the corresponding eigenvector,
choose $\zeta=\zeta_0$ and $\ x=x_0$ in the FMI. Then
$$ D((\zeta_0I-T)x_0,(\zeta_0I-T)x_0)=0 \ , $$
and hence $(s(\zeta_0)E-M)x_0=0$, that is $s(\zeta_0)Ex_0=Mx_0$.
Thus, the value $s(\zeta_0)Ex_0$ is the same for any solution
$s(\zeta)$ of the FMI.\par
To formulate the theorem we need to define the Sz.-Nagy--Foias
function space $\kline_s$ (see \(4\), \(5\)).\par
For $s\in B(\lline,\lline^{\prime})$ the space $\kline_s$ is the
set of all vector-valued functions $f=f_+\oplus f_-$ which
satisfy the following two conditions:
$$ \eqalign{1)\ & f_+(\zeta)\in H^2_+(\lline^{\prime}) {\rm \ \
and\ \ } f_-(\zeta)\in H^2_-(\lline)\ , {\rm \ where\ }
H^2_+(\lline^{\prime}) {\rm \ and\ } H^2_-(\lline)\ {\rm are\
vector\ } \cr
&{\rm Hardy\ spaces\ of\ the\ disc\ } \vert\zeta\vert \leq 1\
{\rm with\ coefficients\ in\ the\ indicated\ space;} \cr
2)\ & {\int_{_{\lt}}\langle\left\(\matrix{I_{\ln^{\prime}}& s\cr
s^* & I_{\ln}\cr}\right\)^{\(-1\)}f,f\rangle dm<\infty ,
{\rm \ where\ } \tline {\rm \ is\ the\ unit\ circle\ and\ } dm\
{\rm\ is}} \cr
&{\rm normalized\ Lebesgue\ measure\ on\ it.} \cr} $$\par
This integral defines an inner product in $\kline_s$. Sometimes
it is convenient to use the following equivalent definition of
the space $\kline_s $ as the set of all vector valued functions
$f=f_+\oplus f_-$  such that
$$ \eqalign{1^{\prime})\ & f_+(\zeta)\ {\rm is\ holomorphic\ in
\ } \vert\zeta\vert<1\ ,\ f_- (\zeta)\ {\rm is\ anti-holomorphic
\ in\ }\vert\zeta\vert<1\ ,\ f_-(0)=0 ,\cr
2^{\prime})\ & \lower 5pt\hbox{$\buildrel \sup\over{0<r<1}$}
\int _ {\lt_r}\langle\left\( \matrix{I_{\ln^{\prime}} & s \cr
s^* & I_{\ln}\cr}\right\)^{\(-1\)}f,f\rangle dm<\infty ,\cr} $$
where $\lt_r$ is the circle of radius $r$ centered at the origin
and $dm$ is normalized Lebesgue measure on it.\par
Note that vector-functions from $\kline_s$ have the following
properties:
$$ \langle\left\(\matrix{I_{\ln^{\prime}}& s \cr s^* & I_{\ln}
\cr}\right\)^{\(-1\)}f,f\rangle(rt) \rightarrow \langle\left
\(\matrix{I_{\ln^{\prime}} & s \cr s^* & I_{\ln} \cr}\right\)^
{\(-1\)}f,f\rangle (t), $$
for almost all $\ t\ ,\ \vert t \vert=1 \ $ (as $r \rightarrow 1$
) and
$$ \int _{\lt_r}  \langle\left\(\matrix{I_{\ln^{\prime}}& s \cr
s^*& I_{\ln}\cr}\right\)^{\(-1\)}f,f\rangle(dm)
\rightarrow \int_{\lt} \langle\left\(\matrix{I_{\ln^{\prime}} & s
 \cr s^* & I_{\ln}\cr}\right\)^{\(-1\)}f,f\rangle dm \ , $$
(as $r \rightarrow 1$).
\noindent
The interpolation sense of  the FMI can be seen from the
following:\par
{THEOREM 1.}~~{\sl Assume that the intersection of the spectrum
of $T$ with the set $\vert\zeta\vert\neq 1$ consists of isolated
points only, that the vector-valued functions $E((\zeta I-T)^{-1}
x)$ and $M((\zeta I-T)^{-1}x)$ are holomorphic for $\vert\zeta
\vert\ne 1$ (where they are defined) and the function
$D(x,\ {I+T\overline\zeta\over I-T\overline\zeta}\ x)$ is
holomorphic in the disc $\vert\zeta\vert<1$ (for all those
$\zeta$ for which it is defined). Then every holomorphic
operator-valued function $s(\zeta)$ which is a solution of the
FMI has the following properties: $F_sx \in\kline_s\ \ $ and
$$ \langle F_s x,F_s x\rangle _{\lk_s}\leq D(x,x), \
(\forall\ x\in\xline)\ , $$
where
$$ \eqalign{& F_sx \ {\buildrel def\over=}\ (F_sx)_+\oplus(F_sx)
_-\in\kline_s, \cr
& (F_sx)_+(\zeta) \ {\buildrel def\over=}\ (s(\zeta)E-M)(\zeta
I-T)^{-1}x, \ \vert\zeta\vert<1 \ , \cr
& (F_sx)_-(\zeta) \ {\buildrel def\over=}\ \overline\zeta
(E-s^*(\zeta)M)(I-T\overline\zeta)^{-1}x , \ \vert\zeta\vert<1
\ . \cr} $$}\par
{PROOF.} \ \ Since $x\in\xline$ and $\ell\in\lline$ are
chosen arbitrarily the  FMI is equivalent to the inequality
$$ \eqalign{D((I-T\overline\zeta)x,(I-T\overline\zeta)x)&+2Re
\langle (E-s^*(\zeta)M)x,\ell\rangle  \cr
&+\lal {I_{\ln}-s^*(\zeta)s(\zeta)\over 1-\overline\zeta\zeta} \
\ell,\ell\rar\geq 0 \ .\cr}\ \eqno(1) $$
We now replace the vector $\ell$ in (1) by the vector
$(\ell-E(I-T\overline\zeta)x)$ and then, after multiplying the
resulting expression through by $1-\overline\zeta\zeta ,$
combine separately the terms which are quadratic with respect
to $x$,  the terms which are linear with respect to $x$ and the
terms which are independent\ of $x$.
\par
The quadratic term has the form
$$ \eqalign{C_2=&(1-\zeta\overline\zeta)D((I-T\overline\zeta)x,
(I-T\overline\zeta)x)-(1-\zeta\overline\zeta)\langle E(I-T
\overline\zeta)x,(E-s^*(\zeta)M)x\rangle\cr\sksix
&~-(1-\zeta\overline\zeta)\rangle (E-s^*(\zeta)M)x,E(I-T
\overline\zeta)x\rangle\cr\skthr
&~+\langle(I_{\ln}-s^*(\zeta)s(\zeta))E(I-T
\overline\zeta)x, E(I-T\overline\zeta)x\rangle  \ . \cr} $$
The linear term is
$$ \eqalign{C_1 & = (1-\zeta\overline\zeta)\langle (E-s^*(\zeta)
M)x,\ell\rangle   \cr
& -\langle (I-s^*(\zeta)s(\zeta))E(I-T\overline\zeta)x,\ell
\rangle  \ . \cr} $$
The constant term is
$$ C_0= \langle (I_{\ln}-s^*(\zeta)s(\zeta))\ell,\ell \rangle
\ . $$
In terms of this notation, inequality (1) can be expressed in
the form
$$ C_2+2ReC_1+C_0\geq 0 \ .  \ \eqno(2) $$
Using the arbitraryness of $x$ and $\ell$, we shall rewrite (2)
in the form
$$ \left\( \matrix{C_2 & \overline C_1 \cr\cr C_1 & C_0 \cr}
\right\)\geq 0 \ . \ \eqno(3) $$
>From the obvious identity
$$ (1-\zeta\overline\zeta) I =(I-T\overline\zeta)-\overline\zeta
(\zeta I-T) , \ \eqno(4) $$
we have
$$ \eqalign{& (1-\zeta\overline\zeta)(E-s^*(\zeta)M)x  \cr\skthr
& = (E-s^*(\zeta)M)(I-T\overline\zeta)x-\overline\zeta(E-s^*
(\zeta)M)(\zeta I-T)x \ . \cr} \eqno(5) $$
It follows from (5) that
$$ \eqalign{(& 1-\zeta\overline\zeta)(E-s^*(\zeta)M)x-(I_{\ln}-
s^*(\zeta)s(\zeta))E(I-T\overline\zeta)x  \cr\skthr
& = s^*(\zeta)(s(\zeta)E-M)(I-T\overline\zeta)x-\overline\zeta
(E-s^*(\zeta)M)(\zeta I-T)x \ ,\cr} \eqno(6) $$
and hence that
$$ C_1=\langle s^*(\zeta)(s(\zeta)E-M)(I-T\overline\zeta)x-
\overline\zeta(E-s^*(\zeta)M)(\zeta I-T)x,\ell \rangle  \ . \
\eqno(7) $$
Next, we transform the second term in $C_2$ with the help of (5)
and the sum of the third and fourth terms with the help of (6) to
obtain:
$$ \eqalignno{C_2 & = (1-\zeta\overline\zeta)D((I-T\overline
\zeta)x,(I-T\overline\zeta)x)  \cr\skthr
&~~-\langle E(I-T\overline\zeta)x,(E-s^*(\zeta)M)(I-T\overline
\zeta)x-\overline\zeta(E-s^*(\zeta)M)(\zeta I-T)x\rangle\cr\skthr
&~~-\langle s^*(\zeta)(s(\zeta)E-M)(I-T\overline\zeta)x-\overline
\zeta(E-s^*(\zeta)M)(\zeta I-T)x,E(I-T\overline\zeta)x\rangle
\cr\sksix
& = (1-\zeta\overline\zeta)D((I-T\overline\zeta)x,(I-T\overline
\zeta)x) & (8) \cr\skthr
&~~+2Re\langle  E(I-T\overline\zeta)x,\overline\zeta(E-s^*(\zeta)
M)(\zeta I-T)x\rangle\cr\skthr
&~~-\langle E(I-T\overline\zeta)x,(E-s^*(\zeta)M)(I-T\overline
\zeta)x\rangle\cr\skthr
&~~-\langle s^*(\zeta)(s(\zeta)E-M)(I-T\overline\zeta)x,E(I-T
\overline\zeta)x\rangle  \ .\cr
\noalign{\noindent{\rm Since}}
\langle s^*(\zeta(s(\zeta&)E-M)(I-T\overline\zeta)x,E(I-T
\overline\zeta)x\rangle\cr\skthr
& =\langle (s(\zeta)E-M)(I-T\overline\zeta)x,s(\zeta)E(I-T
\overline\zeta)x\rangle\cr
\noalign{\noindent{\rm and}}
\langle E(I-T\overline\zeta&)x,(E-s^*(\zeta)M)(I-T\overline\zeta)
x\rangle\cr\sksix
& =\langle E(I-T\overline\zeta)x,E(I-T\overline\zeta)x\rangle
-\langle E(I-T\overline\zeta)x, s^*(\zeta)M(I-T\overline\zeta)x
\rangle\cr\sksix
& =\langle E(I-T\overline\zeta)x,E(I-T\overline\zeta)x\rangle
-\langle s(\zeta)E(I-T\overline\zeta)x,M(I-T\overline\zeta)x
\rangle\cr\sksix
& =\langle E(I-T\overline\zeta)x,E(I-T\overline\zeta)x\rangle
-\langle M(I-T\overline\zeta)x,M(I-T\overline\zeta)x\rangle
\cr\skthr
&~~-\langle (s(\zeta)E-M)(I-T\overline\zeta)x,M(I -T\overline
\zeta)x\rangle  \ , \cr} $$
the expression (8) for $C_2$ takes the form
$$ \eqalign{C_2 & = (1-\zeta\overline\zeta)D((I-T\overline\zeta)
x,(I-T\overline\zeta)x)\cr\skthr
& + 2Re\langle  E(I-T\overline\zeta)x,\overline\zeta(E-s^*(\zeta)
M)(\zeta I-T)x\rangle\cr\skthr
& -\langle E(I-T\overline\zeta)x,E(I-T\overline\zeta)x\rangle+
\langle  M(I-T\overline\zeta)x, M(I-T\overline\zeta)
x\rangle\cr\skthr
& - \langle s(\zeta)E-M)(I-T\overline\zeta)x,(s(\zeta)E-M)(I-T
\overline \zeta)x\rangle  \ . \cr} \eqno(9) $$
Using the Fundamental Identity  \ and grouping similar terms in
the quadratic form $D$ we can reexpress (9) as
$$ \eqalign{C_2 & = ReD((T+\zeta I)(I-T\overline\zeta)x,(T-\zeta
I)(I-T\overline\zeta)x)\cr\skthr
& + 2Re\langle  E(I-T\overline\zeta)x,\overline\zeta(E-s^*(\zeta)
M)(\zeta I-T)x\rangle   \cr\skthr
& - \langle (s(\zeta)E-M)(I-T\overline\zeta)x,(s(\zeta)E-M)
(I-T\overline \zeta)x\rangle  \ . \cr} \eqno(10) $$
Next, upon taking into account expression (7) for $C_1$ and
block-matrix lemma, one can transform inequality (3) to the
following equivalent inequality:
$$ C_2 \geq \Vert (I_{\ln}-s^*(\zeta)s(\zeta))^{\(-1/2\)}
\( s^*(\zeta)(s(\zeta)E-M)(I-T\overline\zeta)x -\overline\zeta
(E-s^*(\zeta)M)(\zeta I-T)x \)\Vert ^2 . \eqno (11) $$
Then, upon substituting formula (10) for $C_2$ into (11) and
using Proposition 3 we obtain
$$ ReD((T+\zeta I)(I-T\overline\zeta)x,\ (T-\zeta I)(I-T\overline
\zeta)x)  $$
$$ + 2Re\langle E(I-T\overline\zeta)x,\overline\zeta(E-s^*(\zeta)
M)(\zeta I-T)x\rangle  \eqno(12) $$
$$ - \lal\left\(\matrix{I_{\ln^{\prime}} & s(\zeta)\cr\cr s^*
(\zeta) & I_{\ln}\cr}\right\)^{\(-1\)}  \left\(\matrix{(s(\zeta)
E-M)(I-T\overline\zeta)x\cr\cr \overline\zeta(E-s^*(\zeta)M)
(\zeta I-T)x \cr}\right\) \ ,
\left\(\matrix{(s(\zeta)E-M)(I-T\overline\zeta)x \cr\cr
\overline\zeta(E-s^*(\zeta)M)(\zeta I-T)x \cr}\right\)\rar\
\geq 0 \ . $$\par\noindent
The inequality (12) can be considered as the final form of the
transformed FMI. We emphasize that all the transformations are
based on identities and do not use any spectral properties of the
operator $T$.\par
The left hand side of the inequality (12) admits a dual
representation. It can be obtained from expression (12) by
regrouping the entries in the first two terms and invoking the
identity
$$\eqalign{{1\over 2}D((T&+\zeta I)y_1,(I-T\overline\zeta)y_2)-
\zeta\langle Ey_1,(E-s^*(\zeta)M)y_2\rangle   \cr\skthr
= & \ {1\over 2}\ D((T-\zeta I)y_1,(I+T\overline\zeta)y_2)+\zeta
\langle (s(\zeta)E-M)y_1,My_2 \rangle \ , \cr} \ \eqno(13) $$
which follows directly from the FI for arbitrary $y_1$ and
$y_2\in\xline$. Inserting $y_1=(I-T\overline\zeta)x$ and $y_2=
(T-\zeta I)y\ $ into (13), we obtain
$$ \eqalign{&{1\over 2}\ D((T+\zeta I)(I-T\overline\zeta)x,(T-
\zeta I)(I-T\overline\zeta)y)\cr\skthr
&~~~+\zeta\langle E(I-T\overline\zeta)x,
\ (E-s^*(\zeta)M)(\zeta I-T)y\rangle\cr\sksix
& = \ {1\over 2}\ D((I-T\overline\zeta)(\zeta I-T)x,(I+T\overline
\zeta )(\zeta I-T)y)\cr\skthr
&~~~-\zeta\langle  (s(\zeta)E-M) (I-T\overline
\zeta)x, M(\zeta I-T)y\rangle \ . \cr} \ \eqno (14) $$
for arbitrary  $x,y\in\xline$ .
Substituting (14) into (12), we obtain
$$ ReD((I-T\overline\zeta)(\zeta I-T)x,(I+T\overline\zeta)(\zeta
I-T)x) $$
$$ -2Re\zeta\langle (s(\zeta)E-M)(I-T\overline\zeta)x,M(\zeta
I-T)x\rangle \eqno (12^{\prime}) $$

$$ -\lal\left\(\matrix{I_{\ln^{\prime}} & s(\zeta)\cr\cr s^*
(\zeta)&I_{\ln}\cr}\right\)^{\(-1\)}  \left\(\matrix{(s(\zeta)E-
M)(I-T\overline\zeta)x\cr\cr \overline\zeta(E-s^*(\zeta)M)(\zeta
I-T)x \cr}\right\) \ , \left\(\matrix{(s(\zeta)E-M)(I-T\overline
\zeta)x\cr\cr\overline\zeta(E-s^*(\zeta)M)(\zeta I-T)x \cr}\right
\)\rar\ \geq 0 \ . $$
\par\noindent
Finally , we shall use the spectral properties of operator $T$ .
Substituting the vector $\ (I-T\overline\zeta)^{-1}(\zeta I-T)^
{-1} x\ $ in place of $x$ in ($12^{\prime}$) we obtain
$$ ReD(x,\ {I+T\overline\zeta\over I-T\overline\zeta}\ x)-
 2Re\zeta\langle (s(\zeta)E-M)(\zeta I-T)^{-1}x,M(I-T\overline
\zeta)^{-1}x\rangle $$
$$ -\ \Vert\left\(\matrix{I_{\ln^{\prime}} & s(\zeta)\cr\cr s^*
(\zeta) &I_{\ln}\cr}\right\)^{\(-1/2\)}  \left\(\matrix{(s(\zeta)
E-M)(\zeta I-T)^{-1}x \cr\cr\overline\zeta(E-s^*(\zeta)M)(I-T
\overline\zeta)^{-1}x \cr}\right\) \Vert ^2\geq 0.\eqno (15) $$
Let us now recall the notation:
$$ \eqalign{& (F_sx)_+(\zeta) \
{\buildrel def\over=}\ (s(\zeta)E-M)(\zeta I-T)^{-1}x\ , \
\cr\skthr
& (F_sx)_- (\zeta)\ {\buildrel def\over=}\
\overline\zeta(E-s^*(\zeta)M)(I-T\overline\zeta)^{-1}x \ ,
\cr\skthr
& F_sx\ {\buildrel def\over=}\ (F_s x)_+\oplus(F_s x)_- \ ,
\cr} $$
and define
$$ P_\zeta(x,y)=\ {1\over 2}\ D(x,\ {I+T\overline\zeta\over I-T
\overline\zeta}\ y) -\zeta \langle (F_s x)_+(\zeta),\ M(I-T
\overline\zeta)^{-1}y\rangle  \ . \ \eqno(16) $$
Then, in view of (14), we also  have
$$ P_\zeta(x,y)=\ {1\over 2}\ D({T+\zeta I\over T-\zeta I}\ x,y)
+\langle E(\zeta I-T)^{-1}x, (F_s y)_-(\zeta)\rangle  \ . \
\eqno(16^{\prime}) $$
\par\noindent
In terms of these notations inequality (15) can be expressed in
the following form:
$$ P_\zeta(x,x)+\overline{P_\zeta(x,x)}- \lal\left\(\matrix{I_
{\ln^{\prime}}& s(\zeta)\cr\cr s^*(\zeta) & I_{\ln} \cr} \right\)
^{\(-1\)}(F_s x)(\zeta), (F_s x)(\zeta)\rar\geq 0 \ . \
\eqno(17) $$
By assumption, the function $P_\zeta(x,x)$ is holomorphic
everywhere  in $\vert\zeta\vert<1$ with the possible exception
of a set of isolated  points. It follows from (17) that the real
part of the function $P_\zeta(x,x)$ is nonnegative. Hence, all
the singularities of this function in the disk $\vert\zeta\vert<
1$ are removable. Furthermore, in view of (17) and the corollary
to Proposition 3, the functions $(F_s x)_\pm(\zeta)$ posess a
harmonic majorant:
$$ P_\zeta(x,x)+\overline{P_\zeta(x,x)}\ge \Vert (F_s x)_\pm(
\zeta) \Vert ^2. $$
This implies that these functions are in $H^2_+(\lline^{\prime})$
\ and\  $H^2_-(\lline)$ respectively. In particular, all their
singularities are removable. Moreover, it can be seen directly
from the definition that $(F_s x)_-(0)=0$.
\par
It follows from formula (16) and from the regularity of the
functions $(F_s x)(\zeta)$ on $\vert\zeta\vert<1$ that
$$ P_0(x,y)=\ {1\over 2}\ D(x,y) \ . $$
Integrating the inequality (17) over the circle $\lt_r$ of
 radius $r$ centered at the origin with respect to normalized Lebesgue measure
$dm(\zeta)={1\over 2\pi i}{d\zeta \over \zeta}$
we get
$$ \eqalign{\int_{\lt_r} \langle\left\(\matrix{I_{\ln^{\prime}} & s(\zeta)
\cr s(\zeta)^* & I_{\ln}\cr}\right\)^{\(-1\)}(F_s x)(\zeta),(F_s x)(\zeta)\rangle dm(\zeta)
& \leq \int_{\lt_r}\(P_{\zeta}(x,x)+\overline{P_{\zeta}(x,x)}\)dm(\zeta) \cr
& = P_0(x,x)+\overline{P_0(x,x)}=D(x,x) \ . \cr} $$
Thus $F_sx\in \kline _s$ and $\langle F_s x,F_s x \rangle _{\kline_s}\ \leq D(x,x)$.
The theorem is proved. \hfill\blacksq
\par
The following proposition shows how the action of the operator $T$
(which was introduced in Theorem 1) is
changed by the  transformation $F_s$ which acts from the space  $\xline$
into the space $\kline_s$.
\medskip
{PROPOSITION 5.} {\sl Let $T$ and $F_s$ be the same as in Theorem 1, then
$$ (F_sTx)(t)\ {\buildrel a.e.\over =}\ t(F_sx)(t)-\left\(\matrix{
I_{\ln^{\prime}}& s(t)\cr\cr s^*(t) & I_{\ln}\cr}\right\)\left\(\matrix{-Mx
\cr\cr Ex}\right\)\ , \ \ (\vert t \vert=1) \ \eqno(18) $$}
\par
{PROOF.} This follows from the definition of $F_s$ by a straightforward
{calculation.\hfill\blacksq}
\par
{REMARK.} \ \ If the spectral condition for $T$ which was formulated
in Theorem 1 is satisfied, then the transformation $F_s$ from $\xline$
in $\kline_s$ is defined uniquely by the relation (18).
\par
The following proposition is a converse of Theorem 1.
\par
{PROPOSITION 6.} {\sl Let $\lline,\ \lline^{\prime},\ \xline,\
T,\ D,\ E\ {\rm and}\ M$ be the objects
occurring}\footnote{$^{(*)}$}{\ninerm We do not impose any spectral
conditions
on $T$ here.} {\sl in Theorem 1, let $s\in B(\lline,\lline^
{\prime})$ and let $(F_s x)(t)
$ be a family of functions in the variable $t \in \tline$ which
depends linearly on $x \in \xline$ and is
defined by the following (generically implicit) formula:
$$ (F_sTx)(t)\ {\buildrel a.e.\over =}\ t(F_s x)(t)-\left\(\matrix{
I_{\ln^{\prime}}& s(t)\cr\cr s^*(t) & I_{\ln}\cr}\right\)\left\( \matrix
{-Mx\cr\cr Ex}\right\),\ (\vert t \vert=1). \leqno {\rm i)} $$
Assume further that
$$\leqalignno{&\ F_s x\in\kline_s \ , \ \forall \ x\in\xline
& {\rm ii)} \cr
\noalign{\noindent and}
& \ \langle F_s x,F_s x\rangle _{\lk_s}\ \leq D(x,x), \ \forall \
x\in\xline \ . & {\rm iii)} \cr }$$
\par\noindent
Then $s(\zeta)$ is a solution of the FMI.}
\par
{PROOF.} Fix a point $\zeta$ with $\vert \zeta \vert < 1$ and consider
the pair of vectors
$$(F_s x)(t) {\rm \ \ and\ \ }
\left\(\matrix{I_{\ln^{\prime}} & s(t) \cr s^*(t) & I_{\ln}\cr}
\right\)\left\(\matrix{I_{\ln^{\prime}}\cr
-s^*(\zeta)\cr}\right\) {\ell^{\prime}\over 1-t\overline\zeta}\ , \ \ell^{\prime}\in\ln^{\prime}\ ,$$
both of which belong to $\kline_s$. Then, upon
calculating all pairwise scalar products (in $\kline_s$) formed
from them and writing out the nonnegativity condition for the Gram matrix we obtain
$$ \left\( \matrix{
\langle F_sx,F_s x\rangle _{\kline_s} & & \overline{\langle (F_sx)_+(\zeta),\ell^{\prime}\rangle } \cr\cr\cr
\langle (F_sx)_+(\zeta),\ell^{\prime}\rangle  && \langle {I_{\ln^{\prime}}-s(\zeta)s^*(\zeta)\over 1-\zeta\overline\zeta}
\ \ell^{\prime},\ell^{\prime}\rangle  \cr} \right\)\geq 0 \ .\ \eqno(19) $$
Substituting the vector $(\zeta I-T)x$ in place of the vector
$x$ in (19) and invoking the analytic continuation of $F_s x$
(which is defined on the boundary in i))
and the linearity of $F_sx$ in $x$, we have
$$ (F_s(\zeta I-T)x)_+(\zeta)=(s(\zeta)E-M)x \ , \ \eqno(20) $$
and, in view of iii),
$$ \langle F_s(\zeta I-T)x,F_s(\zeta I-T)x\rangle _{\lk_s}\ \leq D((\zeta I-T)x,(\zeta I-T)x) \ .
\ \eqno(21) $$
Inserting (20) and (21) into (19) we obtain the FMI$^{\prime}$.
\par
The FMI can be obtained analogously by considering the Gram matrix of the pair of vectors
$$(F_sx)(t) {\rm \ \ and \ \ }
\zeta\left\(\matrix{I_{\ln^{\prime}} & s(t) \cr s(t)^* & I_{\ln}\cr}
\right\)\left\(\matrix{-s(\zeta)\cr I_{\ln}\cr}\right\)
{\ell\over t-\zeta}\ , \ \ell\in\ln \ . \eqno\blacksq $$
\medskip
\noindent{\bf 3. THE ABSTRACT INTERPOLATION PROBLEM}
\medskip
Let $\xline$ be a linear space, $T$ a linear operator on $\xline,\ D$ a
nonnegative sesquilinear form in $\xline$ and let \ $E$ and $M$ be linear operators
from $\xline$ into the Hilbert spaces $\lline$ and $\lline^{\prime},$ respectively.
Suppose, moreover, that the identity
$$ D(x,x)-D(Tx,Tx)=\ \langle Ex,Ex\rangle _{\ln}-\ \langle Mx,Mx\rangle _{\ln^{\prime}} $$
is satisfied.
\medskip
\noindent
{\sl Let $(F_s x)(t)
$ be a family of functions in variable $t \in \tline$ which
depends linearly on $x \in \xline$ and is
defined by the following (generically implicit) formula:
\bigskip\noindent
i)~~$(F_sTx)(t)\ {\buildrel a.e.\over =}\ t(F_s x)(t)-\left\(\matrix{
I_{\ln^{\prime}}& s(t)\cr\cr s^*(t) & I_{\ln}\cr}\right\)\left\( \matrix
{-Mx\cr\cr Ex}\right\),\ (\vert t \vert=1).$\bigskip\noindent
The operator--valued function $s(\zeta) \in B(\lline,\lline^{\prime})$
which is holomorphic in the disc $\vert\zeta\vert<1$ is said to be a solution
of the Abstract Interpolation Problem if
\par\noindent
ii)~~$F_s x\in\kline_s$ \ , \ $\forall \ x\in\xline$
\par\noindent
and\par\noindent
iii)~~$\langle F_s x,F_s x\rangle _{\lk_s}\ \leq D(x,x), \
\forall \ x\in\xline\ .$\par\noindent
Our objective is to describe all the solutions of the Abstract Interpolation Problem.}
\footnote{*}{\ninerm For some choices of data there exists a unique
linear mapping $F_s$ from the space $\xline$
into the space $\kline_s$ with properties i)--iii) for any solution $s(\zeta)$ of the
Abstract Interpolation Problem, for some other data
there might be many mappings
for the same solution $s(\zeta)$. In any case all these mappings $F_s$ can be described along with
the description of the solutions (see [7], [12]).}

The remaining
discussion depends on the paper \(6\). To study the problem
in question, we need some objects connected with the given spaces and
operators. Let $Dx$ denotes the conjugate linear functional defined by the
formula
$$ Dx(y)\ {\buildrel def\over=}\ D(x,y) \ . $$
The scalar product is naturally defined on the set $\{Dx\}_{x\in\lx}$ by the rule
$$ \langle Dx_1,Dx_2\rangle \ {\buildrel def\over=}\ D(x_1,x_2) \ . $$
Obviously this inner product is well-defined. Let us denote by $\kline$
the completion of the space $\{Dx\}_{x\in\lx}$ with respect to the inner
product introduced above.
Then, $\kline$ is a Hilbert space. The Fundamental Identity  \ enables us to
define an isometric operator from the space $\kline\oplus\lline$
into the space $\kline\oplus\lline^{\prime}$. Let us define an operator $V$ by the
formula
$$ V(DTx\oplus Ex)\ {\buildrel def\over =}\ Dx\oplus Mx \ .
\ \eqno(22) $$
The domain $(D_V)$ of the operator $V$ is the closure in
$\kline\oplus\lline$ of the set of all the vectors $DTx\oplus Ex$, the
range $\Delta_V$ is the closure in $\kline\oplus\lline^{\prime}$ of the set of all
the vectors $Dx\oplus Mx$.
\par
Let $\hline$ be a Hilbert space and let $U$ be a unitary operator
from $\hline\oplus\lline$ onto $\hline\oplus\lline^{\prime}$. Following the
paper \(6\) we define the scattering function $s(\zeta)$ of $U$ with
respect to the spaces $\lline$ and $\lline^{\prime}$ in the following way:
$$ s(\zeta)=P_{\ln^{\prime}}U(I_{\lh\oplus\ln}-\zeta P_{\lh}U)^{-1}{\vert\ln} \ .$$
Consider also the functional representation of the space $\hline$
which is defined by formula
$$ Gh=(Gh)_+\oplus(Gh)_- \ , $$
where
$$ \eqalignno{(Gh)_+(\zeta) & =P_{\ln^{\prime}}U(I_{\lh\oplus\ln}-\zeta P_{\lh}
U)^{-1}h  \cr\skthr
\noalign{{\noindent and}}
(Gh)_-(\zeta) & = \overline\zeta P_{\ln}U^*(I_{\lh\oplus\ln^{\prime}}-
\overline\zeta P_{\lh}U^*)^{-1}h
 \ .  & (\vert\zeta\vert<1) \cr} $$
\par
{PROPOSITION 7.} {\it $G$ maps $\hline$ into $\kline_s$,
and}
$$ \Vert Gh\Vert_{\lk_s}\leq\Vert h\Vert_{\lh} \ . $$
\par
The following statement yields a connection between all the
solutions of the abstract interpolation problem and all the
scattering functions of the unitary extensions of the isometry
$V$ (see (22)) with respect to the spaces $\lline$
and $\lline^{\prime}$.
\par
{PROPOSITION 8.} \ \ Let $\hline\supset\kline ,$
let $U$ be unitary
operator from $\hline\oplus\lline$ onto $\hline\oplus\lline^{\prime}$ which extends
$V$ and let $s(\zeta)$ be the scattering function of $U$ with
respect to  the spaces $\lline$ and $\lline^{\prime}$.
Then the functional transformation $F_s$
$$ F_sx\ {\buildrel  def\over=}\ GDx \qquad (x\in X) $$
has property (18), i.e.,
$$ F_sTx= tF_sx - \left\(\matrix{
I_{\ln^{\prime}}& s \cr s^*  & I_{\ln}\cr}\right\)\left\({-Mx
\atop Ex}\right\)\ . $$
\par
{COROLLARY.} {\it The scattering function  of any unitary
extension  of an isometry $V$ is a
solution of the abstract interpolation problem.}
\par
{PROPOSITION 9.} \ \ {\it
Let $s(\zeta)\in B(\lline,\lline^{\prime})$, and let $F_s$ be
a mapping from $\xline$ into $\kline_s$ which satisfies the conditions i)--iii)
of the abstract interpolation problem.
Then $s(\zeta)$ is the scattering matrix of a unitary extension of the
isometry $V$ with respect to the spaces $\lline$ and
$\lline^{\prime}$.}
\par
{COROLLARY.} \ \ {\it The set of all the solutions of the abstract
interpolation problem admits the following description
$$ s(\zeta)=s_{12}(\zeta)+s_{11}(\zeta)\varepsilon(\zeta)\(I-s_{21}(\zeta)\varepsilon(\zeta)\)^{-1}
s_{22}(\zeta),\ \vert\zeta\vert<1 \ , $$
where
$$ S(\zeta)=\left\(\matrix{s_{11}(\zeta)  & s_{12}(\zeta) \cr\cr
s_{21}(\zeta) &  s_{22}(\zeta) \cr} \right\) $$
is {\it the scattering matrix of the isometry}
 $V$, {\rm \(6\)} and $\varepsilon(\zeta)$
is an arbitrary holomorphic contractive operator-valued function which
acts from $\mline_V=(\kline\oplus\lline)\ominus D_V$ into $\nline_V=(
\kline\oplus\lline^{\prime})\ominus\Delta_V$ .}
\bigskip
\noindent{\bf REFERENCES}
\item{1.}Kovalishina I.V. and Potapov V.P., \ \ Indefinite Metric in the
Nevanlinna-Pick Problem. \ \ Dokl. Akad. Nauk Armjan. SSr. {\bf 59},
(1974), no.1, 17-22.
\item{2.}Kheifets A.Ya., Yuditski\v i P.M.,\ \ Interpolation of Operators
commuting with Truncated Shift by Functions of the Schur Class. \ \ Teoriya funkts%
i\v i, funktsional. analys i ikh prilo\^seniya. {\bf 40} (1983), 129-136.
\item{3.}Katsnelson V.E., \ \ Fundamental Matrix Inequality of the Problem
of Decomposition of positive definite Kernel on elementary Kernel. \ \
Kharkov, 1984; \ \ Deposited in UkrNIINTI 10.7.1984, no. 1184, Uk Dep.
\item{4.}Sz.-Nagy B., Foias C., \ \ Analyse Harmonique des Operateurs de
l'Espace de Hilbert. \ \  Ak\'ad\'emiai Kiado, Budapest 1967.
\item{5.}Pavlov B.S., \ \ Selfadjoint Dilatation of a dissipative Operator
and Expansion by its Eigenfunctions. \ \ Matem. Sbornik {\bf 102}, (1977),
no. 4, 511-536.
\item{6.}Arov D.Z., Grossman L.Z., \ \ Scattering Matrix in the Extension
Theory of Isometric Operators. \ \ Dokl. Akad. Nauk. SSSR. {\bf 270},
(1983), No.1, 17-20.
%
%\bye
%  Questions:
%  1.  \lline   or \ln
\medskip
\noindent{\bf LATER REFERENCES ADDED IN TRANSLATION}
\medskip
\item{7.}A.Ya. Kheifets, ``Parseval equality in abstract
interpolation problem and coupling of open systems", {\sl Teor. Funk.,
Funk. Anal. i ikh Prilozhen\/} {\bf 49} (1988) 112--120, {\bf 50}
(1988) 98--103, Russian. English transl., {\sl J. Sov. Math.\/} {\bf
49}, 4 (1990) 1114--1120, {\bf 49}, 6 (1990) 1307--1310.
\medskip
\item{8.}A.Ya. Kheifets, "Generalized bitangential
Schur- Nevanlinna- Pick problem, related Parseval equality and scattering
operator", deposited in VINITI, 11.05.1989, No. 3108--B89 Dep.,
1--60, 1989, Russian.
\medskip
\item{9.}A.Ya. Kheifets, Generalized bitangential
Schur-Nevanlinna-Pick problem and the related Parseval equality, {\sl
Teor. Funk., Funk. Anal. i ikh Prilozhen.\/} {\bf 54} (1990) 89--96,
Russian. English transl., {\sl J. Sov. Math.\/} {\bf 58}, 4 (1992)
358--364.
\medskip
\item{10.}A.Ya. Kheifets, ``Nevanlinna- Adamjan- Arov- Krein
theorem in semi- determinate case", {\sl Teor. Funkt., Funkt. Anal. i
ikh Prilozhen.\/} {\bf 56} (1991) 128--137, Russian. English transl.,
{\sl Journal of Mathematical Sciences\/} {\bf 76}, 4 (1995)
2542--2549.
\medskip
\item{11.}A. Ya. Kheifets ``Scattering Matrices and Parseval Equality
in Abstract Interpolation Problem", Ph. D. thesis, 1990, Kharkov, Russian
\medskip
\item{12.}A. Ya. Kheifets,  P.M. Yuditskii,  ``An analysis  and extension
   of V.P. Potapov's
   approach  to   interpolation   problems with
   applications  to  the generalized  bi-tangential Schur- Nevanlinna- Pick
   problem  and  $j$-inner-outer  factorization", {\sl in Operator Theory:
   Advances and Applications\/}, {\bf 72} (1994) 133--161, Birkhauser
   Verlag, Basel.

\vskip 1truecm\noindent
{\bf Addresses for the authors:}\bigskip
\parindent 0pt\baselineskip 11pt
V. E. Katsnelson\bk Department of Theoretical Mathematics\bk The
Weizmann Institute of Science\bk Rehovot 76100, ISRAEL\smallskip
E-mail: victor.katsnelson@weizmann.ac.il \bigskip
A. Ya Kheifets\bk Department of Mathematics, University of
Massachusetts, Lowell, MA 01854, USA \smallskip E-mail:
 alexander${}_{-} $kheifets@uml.edu
\bigskip

P. M. Yuditskii\bk Department of Mathematics, Bar Ilan University,
52900 Ramat Gan, Israel \smallskip E-mail: yuditski@macs.biu.ac.il
\vskip 1truecm\noindent
AMS Mathematics Subject Classification:\par\noindent
Primary: 47A57, 47A20, 30D50\par\noindent
Secondary: 47A45, 47A48, 30C80
\bye